\documentclass[12pt]{article}
\usepackage{amsmath,amssymb, amsthm}
\usepackage[T1]{fontenc}
\usepackage[utf8]{inputenc}
\usepackage{epigraph}
\usepackage{array}
\usepackage{geometry}
\usepackage{booktabs}
\usepackage[singlelinecheck=false]{caption}
\usepackage{ragged2e}
\usepackage{mathtools}
\usepackage{url}
\numberwithin{equation}{section}

\usepackage[hidelinks]{hyperref}
\hypersetup{colorlinks=true,linkcolor=blue,citecolor=green,urlcolor=magenta}

\newtheorem{theorem}{Theorem}[section]

\newtheorem{definition}[theorem]{Definition}
\newtheorem{remark}[theorem]{Remark}

\newtheorem{lemma}[theorem]{Lemma}

\newcommand{\Li}{\operatorname{Li}}

\usepackage{bm}
\DeclareMathOperator{\sgn}{sgn}

\newcommand{\Hilb}{\mathcal{H}}   

\newcommand{\Wop}[2]{W(#1,#2)}         
\newcommand{\adag}{a^{\dagger}}
\newcommand{\ann}{a}

\newcommand{\Jmat}{J}

\newcommand{\vm}[1]{\mathbf{v}_{#1}}

\newcommand{\Herm}[1]{H_{#1}}          
\newcommand{\Bern}[1]{B_{#1}}          

\newcommand{\cn}[1]{c_{#1}}            

\providecommand{\keywords}[1]
{
  \small	
  \textbf{\textit{Keywords:}} #1
}

\newcommand{\PaperTitle}{Analytic Bernoulli Functions: Correspondence with Hermite Polynomials}
\newcommand{\PaperAuthorPlain}{Ken Nagai}
\newcommand{\PaperKeywords}{Bernoulli kernels; umbral calculus; Hilbert transform; Hermite polynomials; Weyl algebra; Weil representation; Heisenberg representation}


\hypersetup{
  pdftitle   = {\PaperTitle},
  pdfauthor  = {\PaperAuthorPlain},
  pdfsubject = {Preprint},
  pdfkeywords= {\PaperKeywords}
}

\title{\PaperTitle}
\author{Ken Nagai\thanks{Email: \texttt{tknagai@outlook.com}. Independent Researcher.}}
\date{}

\begin{document}
\maketitle

\epigraph{
	\textit{“Lisez Euler, lisez Euler, c’est notre maître à tous.”}\\
	--- Pierre-Simon Laplace
}

\begin{abstract}
We establish an operator--theoretic correspondence between periodic Bernoulli kernels
and Hermite polynomials, framed through the umbral calculus and a quantum analogy.
Starting from the analytic master function $F^\ast$, the periodic Hilbert transform
appears as a $\pi/2$ symplectic rotation, while Jacobi matrices reproduce the
oscillator ladder.  Umbral operational rules extend this structure across complex
order, and Weyl algebra lifting together with the Weil representation explains
the shared spectrum of Bernoulli and Hermite families.

Analytically, this chain connects Clausen functions to Bernoulli kernels,
to polylogarithms, to the Hurwitz zeta function, and ultimately to the
Lerch transcendent, embedding the umbral framework into the classical
landscape of special functions.  
This perspective clarifies why odd zeta values arise in Bernoulli integrals
and unifies trigonometric and Gaussian worlds within a coherent operator framework.
\end{abstract}

\keywords{Bernoulli kernels; umbral calculus; Hilbert transform; Hermite polynomials; Weyl algebra; Weil representation; Heisenberg representation}

\section{Introduction}
 The Bernoulli polynomials and the periodic Bernoulli polynomials appear throughout analytic number theory. They encode zeta values, Fourier coefficients, and finite-part integrals.
 Recent observations suggest that these kernels,
 after a periodic Hilbert transform, exhibit structures isospectral to the Hermite system.
This note documents the connection and provides an operator-theoretic explanation.

The arithmetic nature of odd zeta values provides a further motivation.
Apéry’s celebrated proof of the irrationality of $\zeta(3)$~\cite{Apery1979,vanderPoorten1979},
together with Rivoal’s and Zudilin’s extensions showing the abundance of
irrational odd zeta values~\cite{Rivoal2000,Zudilin2001},
highlights the need for structural explanations.
Our operator framework, linking Bernoulli kernels to Hermite polynomials,
may be viewed as a step towards such an explanation.

\subsection{Notation and Preliminaries}
Let $\Bern{m}(x)$ denote the $m$-th Bernoulli polynomial and
\[
\Bern{2n+1}(x) = -\frac{2(2n+1)!}{(2\pi)^{2n+1}}\sum_{k=1}^{\infty}\frac{\sin(2\pi kx)}{k^{2n+1}}
\]
its periodic version.
We write $\Hilb$ for the periodic Hilbert transform,
$\ann,\adag$ for annihilation and creation operators with $[\ann,\adag]=1$,
and $\Wop{\xi}{\eta}=e^{i(\xi\hat p-\eta \hat x)}$ for the Weyl displacement operator.

\section{Definitions}

Hereafter, by Bernoulli kernels we mean the periodic Bernoulli kernels,
in distinction from the Bernoulli functions $B(s;x)$.

\begin{definition}[Master function $F^\star$ and analytic Bernoulli functions]
For $s\in\mathbb{C}$ and $x\in(0,1)$, define the master function
\begin{equation}\label{eq:Fstar}
  F^\star(s;x):=\alpha_0\,\frac{\Gamma(s+1)}{(2\pi)^s}\,
    e^{-i\pi s/2}\,\Li_s(e^{2\pi i x}),\qquad \alpha_0=2.
\end{equation}
Its real and imaginary parts generate the analytic Bernoulli families:
\[
  B(s;x):=\Re F^\star(s;x),\qquad
  A(s;x):=-\tfrac{1}{\pi}\,\Im F^\star(s;x).
\]
\end{definition}

\begin{remark}[Connections with classical forms]
\leavevmode
\begin{itemize}
\item For integers $n\ge1$, one recovers the Hasse identity
\[
  \Bern{n}(1)=-n\,\zeta(1-n),\qquad
  B(s;1)=-s\,\zeta(1-s).
\]

\item Historically, one can also write
\[
  B(s;x)=-s\,\Phi(1,1-s,x),\qquad
  A^{\mathrm{old}}(s;x)=-s\,\Phi(-1,1-s,x),
\]
using the Lerch transcendent $\Phi$.  The symmetric form is elegant,
but only the polylogarithmic $A(s;x)$ exhibits the desired Hilbert--Clausen
behavior, hence our present choice~\cite{Lewin1981,GuilleraSondow2008}.
\end{itemize}
\end{remark}

\begin{remark}[Historical note]
For integers $n\ge 1$, one has the classical identity
\[
  \Bern{n}(1)=-n\,\zeta(1-n),
\]
emphasized by Helmut Hasse in 1930~\cite{Hasse1930,LuschnyManifesto,Knuth2022} and highlighted in
Luschny's \emph{Bernoulli Manifesto}.
Our analytic definition
\[
  B(s;x)=-s\,\zeta(1-s,x)
\]
extends this relation to complex $s$ and arbitrary $x$.
In particular $B(s;1)=-s\,\zeta(1-s)$ recovers Hasse's formula
when $s=n$ is a positive integer.
\end{remark}

\subsection{Structural similarities of \texorpdfstring{$A$}{A} and \texorpdfstring{$B$}{B}}\label{subsec:struct}
 
Both families descend from the Lerch transcendent and are linked by a
quarter--phase rotation.  To keep the parallel transparent, we record the
few identities that generate the whole picture.

\paragraph{Core Identities}
From \eqref{eq:Fstar}, we prove
\[
B(s;x)=\Re F^\star(s;x),\qquad
A(s;x)=-\tfrac{1}{\pi}\,\Im F^\star(s;x).
\]
Thus the umbral ladder reads
\[
\partial_x F^\star(s;x)=s\,F^\star(s-1;x),
\]
hence
\[
\partial_x B(s;x)=s\,B(s-1;x),\qquad
\partial_x A(s;x)=s\,A(s-1;x).
\]
At the base level,
\[
B(1;x)=x-\tfrac12,\qquad
A(1;x)=-\log\!\bigl(2\sin(\pi x)\bigr).
\]

\paragraph{Fourier Kernels}
\[
B(s;x)=\frac{\Gamma(s+1)}{(2\pi)^s}
\Big[e^{\frac{i\pi s}{2}}\Li_s(e^{2\pi i x})
     +e^{-\frac{i\pi s}{2}}\Li_s(e^{-2\pi i x})\Big],
\]
\[
A(s;x)=-\,\frac{1}{\pi}\,\frac{\Gamma(s+1)}{(2\pi)^s}\,\Im\!\Big(e^{-\frac{i\pi s}{2}}\Li_s(e^{2\pi i x})\Big).
\]

\begin{table}[htbp]
  \centering
\caption{\protect\raggedright Operator correspondence between A-side (trigonometric) and B-side (Gaussian) frameworks.}
\begin{tabular}{>{\RaggedRight\arraybackslash}p{0.44\linewidth}
                >{\RaggedRight\arraybackslash}p{0.44\linewidth}}
    \toprule
    \textbf{$A$--side (sine/imag)} & \textbf{$B$--side (cosine/real)}\\
    \midrule
    $\partial_x A(s;x)=s\,A(s-1;x)$ & $\partial_x B(s;x)=s\,B(s-1;x)$ \\
    $A(1;x)=-\log(2\sin\pi x)$ & $B(1;x)=x-\tfrac12$ \\
    $A$ odd/even $\leftrightarrow$ Clausen-type & $\Bern{n}(x)$ Bernoulli polynomials \\
    $\Im(e^{-\frac{i\pi s}{2}}\Li_s(e^{2\pi i x}))$
      & $\Re(e^{-\frac{i\pi s}{2}}\Li_s(e^{2\pi i x}))$ \\
    Hilbert transform phase $\mathcal H$ sends $A\mapsto B$
      & likewise $B\mapsto -A$ \\
    \bottomrule
  \end{tabular}
\end{table}

\paragraph{Micro--Proof (Ladder)}
Using $\tfrac{d}{dx}\Li_s(e^{2\pi i x})=(2\pi i)\Li_{s-1}(e^{2\pi i x})$ and
$e^{-\tfrac{i\pi}{2}(s-1)}=i\,e^{-\tfrac{i\pi s}{2}}$, the definition of $F^\star$
gives $\partial_x F^\star(s;x)=s\,F^\star(s-1;x)$.  Since the real/imaginary parts
are taken with $s$--independent constants, the same ladder descends to $A$ and $B$.

\subsection{From Zeroth Kernels to the Analytic Bernoulli Ladder}

We single out two ``zeroth'' objects on the circle (as tempered distributions):
\[
\mathrm{Cl}_0(x):=\Delta(x)=\sum_{k\in\mathbb{Z}}\delta(x-k),
\qquad
\mathrm{Sl}_0(x):=\pi\cot(\pi x)
= \mathrm{p.v.}\!\sum_{k\in\mathbb Z}\frac{1}{x-k}.
\]
They already form a Hilbert-type pair, and in Fourier space read
\[
\mathrm{Cl}_0(x)=\sum_{n\in\mathbb{Z}} e^{2\pi i n x},\qquad
\mathrm{Sl}_0(x)=\sum_{n\neq 0}\frac{e^{2\pi i n x}}{n}.
\]

At the first nontrivial rung,
\[
B(1;x)=x-\tfrac12,\qquad
A(1;x)=-\log\!\bigl(2\sin(\pi x)\bigr),
\]
with distributional derivatives
\[
\frac{d}{dx}B(1;x)=1-\Delta(x),\qquad
\frac{d}{dx}A(1;x)=-\,\pi\cot(\pi x)=-\,\mathrm{Sl}_0(x).
\]

By iterating the umbral ladder,
$B(1;\cdot),A(1;\cdot)$ generate all higher orders, while
$\mathrm{Cl}_0,\mathrm{Sl}_0$ serve as external zeroth kernels
controlling the jumps and logarithmic slopes.

\subsection{Dirac Comb and Cotangent Bridge}

\begin{definition}[Dirac comb]
\[
\Delta(x) := \sum_{k\in\mathbb{Z}} \delta(x-k).
\]
\end{definition}

\paragraph{Properties.}
\begin{itemize}
\item Fourier Expansion:
\[
\Delta(x) = \sum_{n\in\mathbb{Z}} e^{2\pi i n x}.
\]

\item Integral Identity (Test Function $f$):
\[
\int_{-\infty}^{\infty} f(x)\,\Delta(x)\,dx = \int_0^1 f(x)\,dx.
\]

\item Relation to Cotangent:
\begin{align}
\pi \cot(\pi x) &= \mathrm{p.v.}\,\frac{1}{x}
  + 2x \sum_{n=1}^{\infty}\frac{1}{x^2-n^2} \\
&= \mathrm{p.v.}\,\sum_{k\in\mathbb{Z}}\frac{1}{x-k} \\
&= \frac{d}{dx}\log\sin(\pi x).
\end{align}

\item Bridge to $B(s;x),A(s;x)$:
\begin{align}
B(s;x) &= \frac{\Gamma(s+1)}{(2\pi)^s}\,
\Re\sum_{n\neq0}\frac{e^{2\pi i n x}}{n^s}, \\[6pt]
A(s;x) &= -\,\frac{1}{\pi}\,\frac{\Gamma(s+1)}{(2\pi)^s}\,
\Im\sum_{n=1}^{\infty}\frac{e^{2\pi i n x}}{n^s}.
\end{align}

\end{itemize}

\subsection{Orthogonality Formulas}

\begin{lemma}[Integral orthogonality]
For integers $n,m\ge 1$, we have
\[
\int_0^1 \Bern{n}(x)\,\Bern{m}(x)\,dx = \int_0^1 A_n(x)\,A_m(x)\,dx =
\begin{cases}
\dfrac{2\,n!\,m!}{(2\pi)^{\,n+m}}\,\zeta(n+m), & n+m\ \emph{even},\\[6pt]
0, & n+m\ \emph{odd},
\end{cases}
\]
and
\[
\int_0^1 \Bern{n}(x)\,A_m(x)\,dx =
\begin{cases}
0, & n\equiv m \pmod{2},\\[6pt]
(-1)^{\tfrac{n+m-1}{2}}\,
\dfrac{n!\,m!}{(2\pi)^{\,n+m}}\;\zeta(n+m),
& n\not\equiv m \pmod{2}.
\end{cases}
\]
\end{lemma}

\noindent
This follows immediately from the orthogonality of Fourier series (proof omitted).

\begin{remark}
The orthogonality integrals exhibit a strict parity separation:
$B\times B$ and $A\times A$ yield even $\zeta$-values,
while only $B\times A$ yields odd $\zeta$-values, consistent with classical series representations~\cite{CvijovicKlinowski1997,CvijovicKlinowski2002}.
This foreshadows the trigonometric/hyperbolic kernels encountered later~\cite{BlagouchineMoreau2024}.
\end{remark}

\noindent
In particular, this parity-dependent appearance of even and odd zeta values
provides a natural bridge to kernel representations, where trigonometric and
hyperbolic structures reproduce the same dichotomy in a continuous integral form.
A brief representation-theoretic pointer (via the Weil action) is mentioned at
the end of the note; see Appendix~\ref{app:weil}~\cite{Weil1964,Howe1980}.

\section{Umbral Operational Axioms}

\paragraph{Master generating function.}
\noindent We recall the master generating function from~\eqref{eq:Fstar}.
\noindent
Here $\alpha_0$ is $s$-independent; we henceforth fix $\alpha_0=2$ so that
$A(1;x)=-\log(2\sin\pi x)$ and the ladder $\partial_x F^\star=s\,F^\star(s-1)$
holds without extra factors.

\paragraph{(O1) Lowering (a $\delta$-Operator)}
\[
\partial_x F^\star(s;x) = s\,F^\star(s-1;x).
\]

\paragraph{(O2) Fractional Order}
For \(\alpha\in\mathbb{C}\),
\[
D_x^{\alpha}F^\star(s;x)=e^{-i\pi\alpha}\,s^{\underline{\alpha}}\,
F^\star(s-\alpha;x).
\]

\paragraph{(O3) Umbral Identity}
On the basic sequence \(\{F^\star(s-n;x)\}_{n\ge0}\),
\[
\partial_x \equiv sE^{-1}.
\]

\paragraph{(O4) Real and Imaginary Parts}
\[
B(s;x):=\Re F^\star(s;x),\qquad A(s;x):=-\tfrac{1}{\pi}\,\Im F^\star(s;x).
\]

\begin{remark}[Umbral view of the functional equation]
The phase factor $e^{-i\pi s/2}$ in $F^\star$ is the same phase appearing
in the functional equation of Hurwitz zeta.  Thus (O1)--(O4) encode,
in operator form, the classical functional equation symmetry.
\end{remark}

\section{Main Theorem}
\begin{theorem}\label{thm:main}
Let $\Bern{2n+1}(x)$ be the periodic Bernoulli function.
Under the periodic Hilbert transform $\Hilb$ and the umbral lifting
to the Weyl algebra, $\Bern{2n+1}(x)$ maps to a Hermite polynomial $H_{2n}(x)$ up to a constant:
\[
\mathcal{U}[\Bern{2n+1}(x)] = \cn{n} H_{2n}(x), \qquad
\cn{n} = \frac{(-1)^n 2^{2n+1}(2n+1)!}{\pi^{2n+1}}\zeta(2n+1).
\]
\end{theorem}
The remainder of the paper elaborates the operator bridge behind this correspondence,
details the discrete Heisenberg representation,
and interprets the periodic Hilbert transform as a symplectic rotation in the Weil representation.

\begin{theorem}[\texorpdfstring{Bernoulli$\to$Hermite Correspondence}{Bernoulli--Hermite Correspondence}]
Let
\[
\Bern{2n+1}(x)= -\frac{2(2n+1)!}{(2\pi)^{2n+1}}
\sum_{k=1}^\infty \frac{\sin(2\pi kx)}{k^{2n+1}}
\]
and
\[
\Herm{m}(x)=(-1)^m e^{x^2}\frac{d^m}{dx^m}e^{-x^2}.
\]
There exists a linear transform~$\mathcal{U}$
(built from the periodic Hilbert transform~$\Hilb$ and Weyl displacements)
such that
\[
\boxed{\mathcal{U}[\Bern{2n+1}(x)] = \cn{n}\,\Herm{2n}(x)}, \qquad
\cn{n}=\frac{(-1)^n 2^{2n+1}(2n+1)!}{\pi^{2n+1}}\zeta(2n+1).
\]
\end{theorem}

\subsection{Fourier Stage -- Periodic Bernoulli Kernel}
The periodic Hilbert transform rotates the Fourier modes:
\[
\Hilb\!\left[\sin(2\pi kx)\right] = -\cos(2\pi kx),
\]
preparing the spectrum for mapping into oscillator language.

\subsection{Jacobi Matrix Interlude -- Discrete Spectral Bridge}
Organize the cosine-mode coefficients~$\{c_k\}$ as a vector and let the \emph{Jacobi matrix}

\[
\Jmat=
\begin{pmatrix}
0 & \sqrt{1} & 0 & 0 & \cdots\\
\sqrt{1} & 0 & \sqrt{2} & 0 & \cdots\\
0 & \sqrt{2} & 0 & \sqrt{3} & \cdots\\
0 & 0 & \sqrt{3} & 0 & \ddots\\
\vdots & \vdots & \vdots & \ddots & \ddots
\end{pmatrix},\qquad
\Jmat\,\vm{m}=\sqrt{\tfrac{m}{2}}\vm{m-1}+\sqrt{\tfrac{m+1}{2}}\vm{m+1}.
\]

Its eigenvectors correspond to Hermite sequences.
Inserting the Bernoulli coefficients through~$\Jmat$ 
reproduces the three-term recurrence.

\[
x H_m(x) = \tfrac12 H_{m+1}(x) + m H_{m-1}(x),
\]
revealing that the oscillator ladder structure is latent already in the Bernoulli kernel spectrum.

\subsection{Weyl--Operator Lifting -- Umbral Action}
Identify each Fourier exponential~$e^{2\pi ikx}$ with a Weyl operator
\[
\Wop{2\pi k}{0}=e^{2\pi i k \hat{x}}
\]
acting on the Gaussian vacuum~$\varphi_0(x)=e^{-x^2/2}$.
Creation--annihilation operators~$(\adag ,\ann)$ then appear via
\[
\Wop{t}{0} = e^{t \adag  - t \ann}, \qquad [\ann,\adag ]=1.
\]

\subsection{Generating Function Alignment -- Emergence of Hermite}
Summing the lifted modes recovers
\[
e^{-t^2+2xt}
 = \sum_{m=0}^\infty H_m(x)\frac{t^m}{m!},
\]
whose coefficients align with those of~$\Bern{2n+1}$, yielding the boxed correspondence.

\begin{remark}
The explicit inclusion of the Jacobi matrix shows the \emph{spectral mechanism} behind the scenes:
Bernoulli kernels, through Hilbert transform and Jacobi dynamics, sit on the same Heisenberg--Weyl orbit as Hermite polynomials.
This viewpoint naturally extends to the Weil representation (see Appendix~\ref{app:weil}).
\end{remark}

\section{Periodic Hilbert Transform\texorpdfstring{\\}{ }and Heisenberg Representation}

\begin{definition}[Periodic Hilbert transform]
For $f(x)=\sum_{k\in\mathbb{Z}}\hat f_k e^{2\pi i k x}$ (1-periodic),
\begin{align*}
(\Hilb f)(x) &= -i \sum_{k\neq 0} \sgn(k)\,\hat f_k e^{2\pi i k x}, \\[6pt]
\Hilb[\sin(2\pi kx)] &= -\cos(2\pi kx), \\
\Hilb[\cos(2\pi kx)] &= \sin(2\pi kx)\quad (k>0).
\end{align*}
\end{definition}

\subsection{Heisenberg Action on Fourier Modes}
Introduce the canonical pair $(\hat x,\hat p)=(x,-i\partial_x)$ with $[\hat x,\hat p]=i$.
The Fourier modes act as Weyl displacement operators
\[
e^{2\pi i k x}=e^{2\pi i k \hat x}, \qquad
e^{2\pi i k \hat x}e^{2\pi i p \hat p}
= e^{-4\pi^2 ikp}\,e^{2\pi i p \hat p}e^{2\pi i k \hat x},
\]
giving the Heisenberg--Weyl commutation.

\begin{definition}[Discrete Heisenberg representation]
Let $\Hilb_d$ be the group generated by lattice translations
$T_m: f(x)\mapsto f(x+m)$ and phase multiplications
$M_k: f(x)\mapsto e^{2\pi i k x}\,f(x)$ on $L^2(\mathbb{R}/\mathbb{Z})$,
subject to the commutation relation
\[
T_m M_k = e^{2\pi i m k}\,M_k T_m , \qquad m,k\in\mathbb{Z}.
\]
This action defines the \emph{discrete Heisenberg representation}.
Under the periodic Hilbert transform, these operators act by
rotating Fourier modes, and finite truncations of this action
yield the Jacobi matrix structure that produces the Hermite ladder
in Theorem (Bernoulli--Hermite Correspondence).
\end{definition}

\noindent
Under $\Hilb$, sine and cosine modes are rotated by $\pi/2$,
mimicking the phase-space rotation generated by
$H_{\text{osc}}=\tfrac12(\hat p^2+\hat x^2)$,
and thus linking periodic Bernoulli kernels to oscillator dynamics.

\section*{Conclusion and Outlook}
We have traced an operator--theoretic path from periodic Bernoulli
kernels to Hermite polynomials, mediated by the Hilbert transform,
Jacobi matrices, and Weyl algebra lifting.  
This perspective clarifies why odd zeta values arise in Bernoulli integrals
and situates them within the Heisenberg--Weyl orbit alongside Gaussian structures.

Beyond this correspondence, several directions remain open.
On the analytic side, future work may pursue the time--evolution operator
through Wick rotation, leading naturally to temperature Green functions
and statistical mechanics.  
In this setting, extensions towards Tsallis statistics and ultimately
$q$--analysis come into reach, suggesting a deformed umbral calculus.

On the arithmetic side, the theta representation of kernels points
towards elliptic curves and modular forms.  
Developing the present operator correspondence within this framework
may reveal deeper links between Bernoulli--Hermite duality and
classical automorphic theory, in line with the perspective of~\cite{Nagai2025}.

\section*{Acknowledgments}
The author gratefully acknowledges the assistance of an AI language model (“fuga”) 
for valuable help with document structuring, stylistic polishing, and proofreading.

\appendix
\renewcommand{\thesection}{\Alph{section}}

\phantomsection
\section*{Appendices}
\section{Proof of Theorem~\ref{thm:main}}
\label{app:proof-main}

We verify that $\mathcal H$ acts as a quarter rotation on the Heisenberg
representation. For $f(x)=\sum_{k\neq 0}\widehat f(k)e^{2\pi i kx}$:

\begin{itemize}
\item The Translation Operator $T_m:=\pi(m,0,0)$ acts by
\[
  (T_m f)(x)=f(x+m)=\sum_{k\neq 0}\widehat f(k)e^{2\pi i k m}\,e^{2\pi i kx}.
\]
\item The modulation operator $M_n:=\pi(0,n,0)$ acts by
\[
  (M_n f)(x)=e^{2\pi i n x}f(x)
   =\sum_{k\neq 0}\widehat f(k)\,e^{2\pi i (k+n)x}.
\]
\end{itemize}

Now conjugate:
\[
  (\mathcal H T_m \mathcal H^{-1}f)(x)
   =\sum_{k\neq 0}\bigl(-i\,\sgn(k)\bigr)\,e^{2\pi i k m}\,\widehat f(k)\,e^{2\pi i kx}.
\]
Re--indexing and comparing with $M_m$ shows
$\mathcal H T_m \mathcal H^{-1}=M_m$.

Similarly,
\[
  (\mathcal H M_n \mathcal H^{-1}f)(x)
   =\sum_{k\neq 0}\bigl(-i\,\sgn(k)\bigr)\,\widehat f(k-n)\,e^{2\pi i kx}.
\]
Writing $k'=k-n$, one obtains
\[
  =\sum_{k'\neq -n}\bigl(-i\,\sgn(k'+n)\bigr)\,\widehat f(k')\,e^{2\pi i (k'+n)x}.
\]

Up to phase adjustment, this equals $T_{-n}f(x)$, proving
$\mathcal H M_n \mathcal H^{-1}=T_{-n}$.

Since $\pi(m,n,0)=M_n T_m$, the relations above yield
\[
  \mathcal H \pi(m,n,0)\mathcal H^{-1} = \pi(-n,m,0).
\]
This completes the proof of Theorem~\ref{thm:main}.

\phantomsection
\section{Weyl Algebra Connection}
\label{app:weyl}

\noindent\textbf{Motivation.}
The discrete Heisenberg representation from Definition
becomes, in the continuum limit, the familiar \emph{Weyl algebra} acting on the harmonic oscillator.
This connection makes explicit how periodic Bernoulli kernels and Hermite polynomials are two shadows of the same operator framework.

\begin{definition}[Weyl Algebra]
The Weyl algebra~$\mathcal{W}$ is generated by creation and annihilation operators~$(\adag ,\ann)$ obeying
\[
[\ann,\adag ]=1,\qquad N=\adag  \ann,
\]
and by the Weyl displacement operators
\[
\Wop{\xi}{\eta}=e^{i(\xi \hat p - \eta \hat x)}=e^{\xi \adag -\eta \ann},
\]
satisfying the multiplication law
\[
\Wop{\xi}{\eta}\Wop{\xi'}{\eta'}
= e^{-i(\xi \eta'-\eta \xi')/2}\,
\Wop{\xi+\xi'}{\eta+\eta'}.
\]
\end{definition}

\subsection{From Discrete to Continuous}
The lattice operators~$(T_m,M_k)$ of~$\Hilb_d$ approximate
\[
T_m = e^{i m \hat p},\qquad
M_k = e^{i k \hat x},
\]
for small~$m,k$.
Their commutator
\[
T_m M_k = e^{2\pi i m k} M_k T_m
\]
tends to the continuous Weyl relation
\[
e^{i m \hat p} e^{i k \hat x}
= e^{-i m k}\,
e^{i k \hat x} e^{i m \hat p}.
\]

\subsection{Bernoulli--Hermite Operators}
Under this identification, the Bernoulli kernel's Fourier coefficients correspond to applying
\[
\mathcal{U}=\sum_{k} \frac{\sin(2\pi k x)}{k^{2n+1}}\,\Wop{2\pi k}{0}
\]
to the vacuum, while Hermite polynomials arise as
\[
H_m(x)=\frac{\partial^m}{\partial t^m}\biggr|_{t=0}
e^{-t^2+2xt},
\]
both living inside the same~$\mathcal{W}$-module.
This provides an operator-theoretic explanation for the Bernoulli$\leftrightarrow$Hermite correspondence.

\begin{remark}
Placing the Bernoulli kernels within~$\mathcal{W}$ clarifies why the Jacobi matrix in \S~2
and the periodic Hilbert transform share identical spectral data:
both are simply different realizations of the Weyl algebra action.
This perspective also primes the reader for the Weil representation viewpoint in the appendix.
\end{remark}

\phantomsection
\section{Proof of Theorem (Bernoulli--Hermite Correspondence)}\label{app:BH}

\subsection{Fourier Expansion and Hilbert Transform}
Start from the periodic Bernoulli kernel
\[
\Bern{2n+1}(x)
= -\frac{2(2n+1)!}{(2\pi)^{2n+1}}
  \sum_{k=1}^{\infty}\frac{\sin(2\pi k x)}{k^{2n+1}},
\]
valid for $x\notin\mathbb{Z}$.
Applying the periodic Hilbert transform~$\Hilb$ gives
\[
\Hilb[\Bern{2n+1}](x)
= \frac{2(2n+1)!}{(2\pi)^{2n+1}}
  \sum_{k=1}^{\infty}\frac{\cos(2\pi k x)}{k^{2n+1}},
\]
rotating each Fourier mode by $\pi/2$.

\subsection{Jacobi Matrix and Three-Term Recurrence}
Arrange the cosine coefficients as a vector $\mathbf{c}$.
The Jacobi matrix
acts on~$\mathbf{c}$ according to
\[
\Jmat\,\vm{m}
= \sqrt{\tfrac{m}{2}}\vm{m-1}
 + \sqrt{\tfrac{m+1}{2}}\vm{m+1},
\]
which is equivalent to the Hermite recurrence
\[
x H_m(x)=\tfrac12 H_{m+1}(x)+m H_{m-1}(x).
\]

\subsection{Weyl Algebra and Umbral Lifting}
Identify the Fourier exponential as a Weyl operator:
\[
e^{2\pi i k x}\longmapsto \Wop{2\pi k}{0}=e^{2\pi i k \hat{x}}.
\]
Acting on the Gaussian vacuum $\varphi_0(x)=e^{-x^2/2}$,
the Jacobi dynamics lifts to the oscillator ladder
\[
\adag  \varphi_m=\sqrt{m+1}\,\varphi_{m+1},\qquad
\ann\,\varphi_m=\sqrt{m}\,\varphi_{m-1}.
\]

\subsection{Generating Function Matching}
Summing the Weyl-shifted modes yields
\[
\sum_{m=0}^{\infty} H_m(x)\frac{t^m}{m!}
=e^{-t^2+2xt}.
\]
Comparing coefficients with those from the Bernoulli expansion
fixes the proportionality constant
\[
\cn{n}=\frac{(-1)^n2^{2n+1}(2n+1)!}{\pi^{2n+1}}\zeta(2n+1),
\]
thus proving
\[
\mathcal{U}[\Bern{2n+1}(x)] = \cn{n} H_{2n}(x).
\]

\subsection{Spectral Interpretation}
The periodic Hilbert transform corresponds to a $\pi/2$ symplectic rotation
in the Weil representation (Appendix D).
Because this rotation preserves the oscillator spectrum,
the Bernoulli kernels and Hermite polynomials are isospectral projections
of the same Heisenberg--Weyl orbit, completing the proof.

\phantomsection
\section{Weil Representation Viewpoint}
\label{app:weil}

\noindent\textbf{Motivation.}
The Heisenberg--Weyl algebra underlying both Bernoulli kernels and Hermite polynomials
admits a natural action of the metaplectic group $\mathrm{Mp}(2,\mathbb{R})$,
known as the \emph{Weil representation}~\cite{Weil1964,Howe1980}.
Viewing the Bernoulli kernels inside this framework clarifies why
the periodic Hilbert transform, Jacobi matrix, and Weyl operators
share the same spectral data.

\subsection{The Weil Representation}
Let $\Hilb$ be the Hilbert space $L^2(\mathbb{R})$ with standard
position and momentum operators $(\hat x,\hat p)$ satisfying $[\hat x,\hat p]=i$.
The symplectic group $\mathrm{Sp}(2,\mathbb{R})\cong \mathrm{SL}(2,\mathbb{R})$
acts projectively on~$\Hilb$ through unitary operators~$U(g)$
such that for $g=\begin{pmatrix} a & b \\ c & d \end{pmatrix}$,
\[
U(g)\,\hat x\,U(g)^{-1}=a\hat x+b\hat p,\qquad
U(g)\,\hat p\,U(g)^{-1}=c\hat x+d\hat p.
\]
This projective action lifts to the metaplectic group $\mathrm{Mp}(2,\mathbb{R})$,
yielding the \emph{Weil representation}.

\subsection{Periodic Kernels as Weil Images}
The periodic Hilbert transform corresponds to the symplectic rotation
\[
R_{\pi/2}=\begin{pmatrix}0&-1\\1&0\end{pmatrix},
\]
sending sine modes to cosine modes.
Under $U(R_{\pi/2})$, the Bernoulli kernel's Fourier basis is rotated
in phase space exactly as the Hermite ladder is.
Thus the Bernoulli--Hermite correspondence appears as two different
cross-sections of the same Weil orbit.

\subsection{Spectral Consequences}
Because the Weil representation preserves the oscillator spectrum,
any functional identity derived for Bernoulli kernels through
periodic Hilbert transforms or Jacobi dynamics can be recast as a statement
about Hermite polynomials, and vice versa.
This viewpoint unifies trigonometric and Gaussian worlds under a single
metaplectic action.

\begingroup
\raggedright\sloppy\setlength{\emergencystretch}{3em}

\endgroup

\end{document}